\documentclass[12pt]{article}
\usepackage{latexsym,amsfonts,amssymb}

\usepackage{amsmath}
\usepackage{amsthm}
\usepackage{amssymb}
\usepackage{delarray}
\usepackage[dvips]{graphics}
\usepackage{epsfig}
\usepackage{color}
\usepackage{setspace}

\newtheorem{thm}{Theorem}[section]
\newtheorem{lemma}[thm]{Lemma}

\theoremstyle{definition}

\def\cd{\cdots}

\def\dim{\text{\rm dim}\,  }
\def\ch{\text{\rm ch}\,  }

\def\hom{\text{\rm Hom}\,  }

\def\al{\alpha}

\def\la{\lambda}

\def\g{\mathfrak g}

\begin{document}
\title{Rigidity of formal characters of  \\
        Lie algebras (II)
    }
\date{}
\author{Zhongguo Zhou  %\quad Jiachen Ye
}
\date{\small\it College of Science, Hohai University\\
    Nanjing,  210098, P.R.China.\\ {
{\rm e-mail: zhgzhou@hhu.edu.cn}\\}
%Department of Applied Mathematics, Tongji University,\\
%Shanghai 200092, People's Republic of China\\
%{\rm e-Mail: jcye@tongji.edu.cn}\\
}
\maketitle
\begin{abstract}
For a complex simple Lie algebra of type $A_l,B_l,C_l$ or $D_l,$
given a family of elements $f_\la\in \mathbb Z[\Lambda],\la\in \Lambda^+,$
we show that $f_\la$ is just the formal character of the
Weyl module $V(\la) $ if $f_\la$ satisfy several natural conditions.
Hence we give a necessary and sufficient condition for constructing
a family of $\g_l$-modules from a family of $\g_{l-1}$-modules. \\


\medskip
\noindent {\bf 2000 MR Subject Classification:} 17B10 \quad 20G05
\end{abstract}
%%%%%%%%%%%%%%%%%%%%%%%%%%%%%%%%%%%%%%%%%%
%%%%%%%%%%%%%%%% section 1 %%%%%%%%%%%%%%%%%%%%%%%%%%
\section{Introduction}

Let $\g_l$ be a complex simple Lie algebra and  $V(\la)$ be the
Weyl module. The formal character of $V(\la)$
$$
\ch_\lambda=\sum_{x\in\Pi^(\lambda)}m_\lambda(x)e(x)
=\sum_{\mu\in\Pi^+(\lambda)}m_\lambda(\mu)h(\mu), \;\;\lambda\in \Lambda^+
$$
is determined by the Weyl character formula or by other methods,
such as the Freudenthal formula or Kostant formula. All the formal
characters $\ch_\lambda, \;\lambda\in \Lambda^+$ is a basis of group
ring $\mathbb Z[\Lambda]^W$ which is invariant under action of Weyl
group $W.$  The Weyl module $V(\la)$ is also a $\g_{l-1}$-module for
the natural Lie subalgebra  $\g_{l-1}$ of Lie algebra $\g_l,\,l\geq 2$.
Especially, we have
%\begin{eqnarray}
$$
\mbox{If}\  \la-\mu=\sum_{i=1}^{l-1}k_i\alpha_i+0\alpha_l, \mbox{then}\
   m_\lambda(\mu)=m_{\lambda|_{\g_{l-1}}}(\mu|_{\g_{l-1}}).\eqno{(*)}
$$
%\end{eqnarray}

The product of $\ch_\mu$ and $\ch_\nu$ is  defined by tensor product
$V(\mu)\bigotimes V(\nu)$ as follow:
$$
\ch_\mu\ch_\nu=\sum_{\la\in \Lambda^+}c_{\mu,\,\nu}^\la\ch_\la.
$$
The Littlewood-Richardson coefficient $c_{\mu,\,\nu}^\la$ defined the
multiplicity of $V(\la)$ in $V(\mu)\bigotimes V(\nu)$ can be
determined according to the formal characters $\ch_\lambda, \ch_\mu$ and $\ch_\nu.$
According to the complete
reducibility of $\g_l$-modules and
$$
\dim\hom(U\otimes V, W)=\dim\hom(U ,W\otimes V^*),
$$
then
%\begin{eqnarray}
$$
c_{\mu,\,\nu}^\la=c_{\la,\,-w_0\nu}^\mu.\eqno{(**)}
$$
%\end{eqnarray}

In \cite{z} we prove a theorem for Lie  algebra $\g_l$ of type $A$ that

\medskip
\textit{ Given a family of elements  $f_\la$'s of  $ \mathbb
Z[\Lambda]^W,\la\in \Lambda^+,$ if the condition (*)and (**) are
satisfied, then these $f_\la$'s  are just equal to the formal characters $\ch_\la$'s.}
\medskip

This theorem describe the relation between $V(\la)$ as $\g_l$-module and as $\g_{l-1}$-module.
That is to say, the first condition (*) tell us the $f_\la$ describe the $\g_{l-1}$-module
structure locally, however the second condition (**) ensure us to construct a
$\g_{l}$-module globally from the local $\g_{l-1}$-module structure.
So the theorem give a  necessary and sufficient condition for lifting
all these $\g_{l-1}$-modules to   $\g_{l}$-modules.
We are especially interested in  finding out a similar condition in positive
characteristic case. This is also the motivation for the paper \cite{z} and this note.

The formal characters $f_\la$'s are determined by these natural condition completely.
This property is called rigidity of formal characters in \cite{z}.
 We continue our work in  \cite{z}  and state a similar rigidity theorem \ref{th1}
 for Lie algebra $\g_l$ of type $A_l,B_l,C_l,\, l\geq 2,$
or $D_l,\, l\geq 4$   in this note with some additional condition.

\section{notations}

From now on let $\g_l$ be a complex simple Lie algebra of type $A_l,B_l,C_l,\, l\geq 2,$
or $D_l,\, l\geq 4,$
and let $\Delta=\{\al_1, \al_2,\cd, \al_l\}$ be the set of simple roots.
Let $\Lambda$ be the set of weights, and $\omega_1, \omega_2,
\cdots, \omega_l$ the set of fundamental dominant weights. Then the
set of dominant weights is denoted by $\Lambda^+.$ Let $\ch_\lambda$,
$\lambda\in \Lambda^+$, be the formal character of the Weyl module
$V(\lambda)$, they form a free $\mathbb Z$-module of the commutative
ring $\mathbb Z[\Lambda]$, with base $\{e(\lambda),\lambda\in\Lambda\},$
and multiplication $e(\lambda)*e(\mu)=e(\lambda+\mu).$  Let $W$ be the
Weyl group, action on  $\mathbb Z[\Lambda]$ naturally as $\sigma e(\lambda)
=e(\sigma\lambda).$ Set
$$
\mathbb Z[\Lambda]^W=\{f\in \mathbb Z[\Lambda]\,|\, wf=f, \;\;w\in W\}.
$$
Let $W_\lambda $ be the $W$-orbit of $\lambda$ and $h(\lambda )=\sum_{x\in
W_\lambda} e(x).$ Let $\Pi(\lambda)$ be the set of saturated weights of
weight $\lambda$ and $\Pi^+(\lambda)=\Pi(\lambda)\cap \Lambda^+.$
It is well known that
$$ \ch_\lambda=\sum_{x\in\Pi^(\lambda)}m_\lambda(x)e(x)=\sum_{\mu\in
\Pi^+(\lambda)}m_\lambda(\mu)h(\mu), \;\;\lambda\in \Lambda^+,$$
forms a basis of $\mathbb Z[\Lambda]^W.$

Recall that $w_0$ is the longest element  of $W$ and $c_{\mu,\,\nu}^\la$  is the
Littlewood-Richardson coefficient.  According to the complete
reducibility of $\g_l$-modules and
$$
\dim\hom(U\otimes V, W)=\dim\hom(U ,W\otimes V^*),
$$
we have
$$
[V(\mu)\otimes V(\la): V(\nu)]=[V(\nu)\otimes V(\la)^*: V(\mu)]=[V(\nu)\otimes V(-w_0\la): V(\mu)].
$$
Hence
$$
c_{\mu,\,\la}^\nu=c_{\nu,\,-w_0\la}^\mu.
$$

\section{main results}

{\bf 3.1.} Let $\beta=\sum_{i=1}^{l}k_i\alpha_i,$ define
$$\mbox{Supp}(\beta)=\{\alpha_i\,| \,k_i>0\}.$$
For $\lambda \in \Lambda^+,$ set
$$
f_{\lambda}
=\sum_{x\in \Pi(\lambda)}n_\lambda(x)e(x)
=\sum_{\mu\in \Pi^+(\lambda)}n_\lambda(\mu)h(\mu)
\in \mathbb Z[\Lambda]^W
$$
  satisfied $n_\lambda(\lambda)=1, n_\lambda(x)=0,
$ if
$x\notin \Pi(\lambda).$ Then $f_{\lambda}$ is also a basis of $\mathbb Z[\Lambda]^W.$
Hence there exists unique $n_{\mu,\,\nu,}^\lambda $ such that
$$
f_{\mu}* f_{\nu}=\sum_{\lambda\in\Pi^+(\mu+\nu)}n_{\mu,\,\nu}^\lambda f_{\lambda}.
$$
By the definition of $f_{\lambda}'s,$ we have
$$n_{\mu,\,\nu}^{\mu+\nu}=1; n_{\mu,\,\nu}^\lambda=n_{\nu,\,\mu}^\lambda.$$

For dominant weights $\lambda,t,\mu,\,\nu, $  the two numbers
$n_\lambda(t) $ and $n_{\mu,\,\nu}^t$ is determined by each other
in some sense as follow c.f \cite{z} 3.2.

\begin{lemma}\label{num}
 If $ \lambda=\mu+\nu, \mu\neq 0,\nu\neq 0$ then
\begin{eqnarray*}\label{eqn3}
n_\lambda(t)=n_{\mu,\,\nu}^t +g(n_\mu(y),n_\nu(z),n_s(x)),
\end{eqnarray*}
where function $g(n_\mu(y),n_\nu(z),n_s(x)) $ is determined by those $ n_\mu(y),
n_\nu(z), n_s(x),$ with $ \mu \lvertneqq \lambda, \nu \lvertneqq \lambda,s\preceq \lambda$ and $
 \mu-y\preceq \lambda-t, \nu-z\preceq\lambda-t,
  s-x\precneqq\lambda-t.$
\end{lemma}

We list some facts on fundamental dominant weights $\omega_i$ c.f.\cite{hum}.

(1)\ For type $A_l,$
\begin{eqnarray*}
&&\omega_1-w_0\omega_1=\omega_1+\omega_l=\omega_l-w_0\omega_l=\alpha_1+\alpha_2+\cdots +\alpha_l;\\
&&\omega_i-w_0\omega_i=\omega_i+\omega_{l-i+1}=\alpha_1+2\alpha_2+\sum_{j=3}^l k_j\alpha_j, \ 1< i<l.\\
\end{eqnarray*}

(2)\ For type $B_l,$
\begin{eqnarray*}
&&\omega_1-w_0\omega_1=2\omega_1=2(\alpha_1+\alpha_2+\cdots +\alpha_l);\\
&&\omega_l-w_0\omega_l=2\omega_l=\alpha_1+2\alpha_2+\cdots +l\alpha_l;\\
&&\omega_i-w_0\omega_i=2\omega_i=2(\alpha_1+2\alpha_2)+\sum_{j=3}^l k_j\alpha_j, \ 1< i<l.\\
\end{eqnarray*}

(3)\ For type $C_l,$
\begin{eqnarray*}
&&\omega_1-w_0\omega_1=2\omega_1=2(\alpha_1+\alpha_2+\cdots +\alpha_{l-1})+\alpha_l;\\
&&\omega_i-w_0\omega_i=2\omega_1=2(\alpha_1+2\alpha_2)+\displaystyle\sum_{j=3}^l k_j\alpha_j, \ 1< i.\\
\end{eqnarray*}

(4)\ For type $D_l,$
\begin{eqnarray*}
&&\omega_1-w_0\omega_1=2\omega_1=2(\alpha_1+\alpha_2+\cdots +\alpha_{l-2})+\alpha_{l-1}+\alpha_{l};\\
&&\omega_l-w_0\omega_l=\omega_{l-1}+\omega_l=\alpha_1+2\alpha_2+\cdots +(l-2)\alpha_{l-2}+(l-1)(\alpha_{l-1}+\alpha_l);\\
&&\omega_i-w_0\omega_i=\alpha_1+2\alpha_2+\displaystyle\sum_{j=3}^l k_j\alpha_j, \ 1< i<l.\\
\end{eqnarray*}
From the definition of ${\rm{Supp}(\beta })$ and the facts on $\omega_i-w_0\omega_i$ in above,
we  have
\begin{lemma}\label{supp}
Let $\g_l$ be a complex simple Lie algebra of type $A_l,B_l,C_l$
or $D_l$ and $\beta=\alpha_1+\alpha_2+\cdots+\alpha_l+\displaystyle\sum_{j=1}^l k_j\alpha_j.$
Then the following statements hold:
\begin{enumerate}
  \item  $|{\rm{Supp}}(\omega_i-w_0\omega_i-\beta)|<l, $  \quad if $k_1\geq 1$ or $k_2\geq 2.$
  \item  $|{\rm{Supp}}(\omega_1-w_0\omega_1-\beta)|<l, $  \quad for $A_l,C_l,D_l.$
  \item $|{\rm{Supp}}(\omega_l-w_0\omega_l-\beta)|<l, $  \quad for $A_l,B_l,D_l.$
  \item $|{\rm{Supp}}(\omega_1-w_0\omega_1-\beta)|<l, $  \quad if $\displaystyle\sum_{j=1}^l k_j>0$ for $B_l.$
  \item $|{\rm{Supp}}(\omega_i-w_0\omega_i-\beta)|<l, $  \quad if $\omega_i$ is the minimal weight.
\end{enumerate}

\end{lemma}
\medskip
{\bf 3.2.}   Let $\lambda,\mu\in \Lambda^+,$ define $\mu< \lambda$ if $\mu\prec \lambda $  or
$ \lambda-\mu\in   \Lambda^+.$

Now we state the main theorem in this paper.
\begin{thm}\label{th1}
Let $\g_l$ be a complex simple Lie algebra of type $A_l,B_l,C_l,\, l\geq 2,$
or $D_l,\, l\geq 4,$ if these $f_{\lambda}$'s satisfy the following  conditions:

(1) $n_\lambda(\mu)=m_\lambda(\mu),$ if $|{\rm{Supp}}(\lambda-\mu)|<l.$

(2) $n_\lambda(\mu)=m_\lambda(\mu),$ if
                  $\la_1=0,$ $\displaystyle\lambda-\mu=\alpha_1+2\alpha_2+\sum_{i=3}^lt_i\alpha_i,t_i\geq 1.$

(3) $n_\lambda(\mu)=m_\lambda(\mu),$ if
                  $\la_1\neq 0,$ $\displaystyle\lambda-\mu=\alpha_1+\alpha_2+\cdots+\alpha_l$ for type $B_l.$

(4) $n_{\mu,\,\nu}^\lambda=n_{\lambda,\,-w_0\nu}^\mu $ for $\la, \mu,\nu\in\Lambda^+.$

\noindent Then $f_{\lambda}=ch_\lambda, n_{\mu,\,\nu}^\lambda=c_{\mu,\,\nu}^\lambda.$
\end{thm}
\proof It is only need to prove $n_\lambda(\mu)=m_\lambda(\mu),
\lambda\in \Lambda^+,\mu\in\Pi^+(\lambda).$ We will prove the
theorem by induction on $\Lambda^+$ with the partial order
$\lq\lq<"$ and on $\Pi(\lambda)$ with the partial order
$\lq\lq\prec".$

Firstly,  if $\lambda=\omega_i$ or  $0$  is a minimal weight,
then their saturated weight set $\Pi(\lambda)$ only contains one
dominate weight. So the theorem holds by the definition of
$f_\lambda.$

Suppose that $\lambda\in \Lambda^+ $ not be a minimal weight. Let
$\mu\in \Pi^+(\lambda),\beta=\lambda-\mu.$ We will consider the different cases
as follows:

(1)\,  When $|\rm{Supp}(\beta)|<l,$ by the first condition in theorem then
$$n_\lambda(\mu)=m_\lambda(\mu).$$

(2)\,  When $|\rm{Supp}(\beta)|=l $ and
$\beta=\alpha_1+\alpha_2+\cdots+\alpha_l+\displaystyle\sum_{j=1}^l k_j\alpha_j.$

Because $\lambda\neq 0$ , there exists $i$ such that $\la_i\neq 0.$ So $\la-\omega_i\in \Lambda^+.$
 By the fourth condition then
$$
n_{\la-\omega_i,\,\omega_i}^\mu
=n_{\mu,\,-w_0\omega_i}^{\lambda-\omega_i}.
$$
Noticing that
$$
\mu+(-w_0\omega_i)-(\lambda-\omega_i)=\omega_i-w_0\omega_i-\beta.
$$

\quad(i)\  When  $k_1\geq 1$ or $k_2\geq 2,$ by lemma \ref{supp},
$|{\rm{Supp}}(\omega_i-w_0\omega_i-\beta)|<l. $ Then by the first condition in the theorem
$$n_{\mu+(-w_0\omega_i)}({\lambda-\omega_i})=m_{\mu+(-w_0\omega_i)}({\lambda-\omega_i}).$$
Moreover  we have
\begin{eqnarray*}
&n_{\mu}(y)=m_{\mu}(y), &  \text{if}\  \mu-y\preceq \omega_i-w_0\omega_i-\beta;\\
&n_{-w_0\omega_i}(z)=m_{-w_0\omega_i}(z), &  \text{if}\  -w_0\omega_i-z\preceq \omega_i-w_0\omega_i-\beta;\\
&n_s(x)=m_s(x) &             \text{if}\  s-x\preceq \omega_i-w_0\omega_i-\beta.
\end{eqnarray*}
Then by lemma \ref{num}
\begin{eqnarray*}
n_{\mu,\,-w_0\omega_i}^{\lambda-\omega_i}&=&n_{\mu+(-w_0\omega_i)}({\lambda-\omega_i})
               -g(n_{\mu}(y),n_{-w_0\omega_i}(z),n_s(x))\\
&=&m_{\mu+(-w_0\omega_i)}({\lambda-\omega_i})
               -g(m_{\mu}(y),m_{-w_0\omega_i}(z),m_s(x))\\
               &=&c_{\mu,\,-w_0\omega_i}^{\lambda-\omega_i}.
\end{eqnarray*}

Hence by the fourth condition in the theorem and the property of $c_{\mu,\nu}^\lambda$ in (**)
$$
n_{\la-\omega_i,\,\omega_i}^\mu=
n_{\mu,\,-w_0\omega_i}^{\lambda-\omega_i}=c_{\mu,\,-w_0\omega_i}^{\lambda-\omega_i}
=c_{\la-\omega_i,\,\omega_i}^\mu.
$$
By lemma \ref{num} and induction hypothesis again, we have
$$n_\lambda(\mu)=m_\lambda(\mu).$$

\quad(ii)\ When $k_1=0$ and $k_2=1$ or $k_1=0$ and $k_2=0$

  If $\la_1=0,$ by the second condition in theorem then
  $$n_\lambda(\mu)=m_\lambda(\mu)$$
  in the first case. However we have $\la_1\neq 0$ in the second
  case because
  $$
\lambda=\mu+\alpha_1+\alpha_2+\sum_{i=3}^lk_i\alpha_i=\mu+(1,\cdots).
  $$.

  If $\la_1\neq 0,$  choose $\lambda_i=\lambda_1,$ then
$$
n_{\la-\omega_1,\,\omega_1}^\mu
=n_{\mu,\,-w_0\omega_1}^{\lambda-\omega_1}.
$$
For
$$
\mu+(-w_0\omega_1)-(\lambda-\omega_1)=\omega_1-w_0\omega_1-\beta,
$$
by lemma \ref{supp}, $|{\rm{Supp}}(\omega_1-w_0\omega_1-\beta)|<l $
for Lie algebra $\mathfrak g_l$ of type $A_l,C_l,D_l;$   type $B_l$
in the first case, and type $B_l$ in the second case with additional
 condition $ \sum_{j=1}^l k_j\alpha_j >0,$ so
$$n_\lambda(\mu)=m_\lambda(\mu)$$
by lemma \ref{num} and induction hypothesis as we prove  in the case (i).

The above equation also holds for type $B_l$ in the second case  when $ \sum_{j=1}^l k_j\alpha_j =0 $
by the third condition in the theorem.

We complete the proof of theorem for Lie algebra $\mathfrak g_l$
of type $A_l,B_l,C_l,D_l.$

{\bf 3.3 Remarks.}
{ 1.}
   As we mentioned in the introduction this theorem describe when we can construct a family of
$\g_l$-modules from a family of $\g_{l-1}$-modules. However it need more conditions for
type $B_l,C_l,D_l$  than type $A_l$ in\cite{z}. This is because there exist some fundamental dominant weights
which are no longer minimal weight.
It need also these conditions  to   obtain a similar theorem for Lie algebra of type $E_l,F_4$ and $G_2.$

{ 2.}
    The number $m_\lambda(\mu)$ in the second and third condition  can be determined precisely, c.f.
 theorem 4 in \cite{wang}, theorem 3.10 in \cite{jyz}  and theorem 4.9 in \cite{yz}.

{ 3.}
  The first condition is very natural and the
  fourth condition can be weaken. In our proof it need only
  $n_{\mu,\,\nu}^\lambda=n_{\lambda,\,-w_0\nu}^\mu $ for $\la, \mu\in\Lambda^+$ and
  $\nu=\omega_i$  be a fundamental dominant weight. We will investigate and
  generalize the weaker condition
 in positive characteristic in the future.

\medskip

%\begin{flushleft}
\noindent{\bf Acknowledgement:} This work was
supported  by the Fundamental Research Funds for the Central
Universities 2009B26914 and 2010B09714.
%\end{flushleft}

%%%%%%%%%%%%%  bibliography      %%%%%%%%%%%%%%%%%%%%%%%%%%%%%
%%%%%%%%%%%%%  bibliography      %%%%%%%%%%%%%%%%%%%%%%%%%%%%%
\makeatletter
\def\@biblabel#1{#1.}
\makeatother


\begin{thebibliography}{9}

\bibitem{gv} Gasharov, Vesselin ,A short proof of the Littlewood-Richardson rule, European Journal of Combinatorics
(1998)19 (4): 451-453,

\bibitem{hum} {J.E.Humphreys,  }  {Introduction to Lie
Algebras and Representation Theory,  {\rm GTM 9}},
{Springer-Verlag},   {New York/Heidelberg/Berlin},
{1972}.

\bibitem{jc}  J. C. Jantzen,  Representations of algebraic groups,
         Academic Press,  Orlando,  1987; 2nded., Amer. Math. Soc., Providence, RI, 2003.

\bibitem{wang}  Wang Yujin,  A  new multiplicity formula for the
Weyl modules of type D, Journal of Hebei University(Natural Science
Edition), (2006)26 (6) 570-573.



\bibitem{jyz} Ye Jiachen, Yun Bensheng, Zhan Jia, A  new multiplicity formula for the
Weyl modules of type B and C, Journal of Mathematical Research $\&$
Exposition  (2010)30(4):675-686.

\bibitem{yz} Ye  Jiachen,  Zhou Zhongguo,  A new multiplicity formula for the
Weyl modules of type $A $  Commun. Algebra, (2005)33(12):4361-4373.


\bibitem{z}   Zhou Zhongguo,  Rigidity of formal character  of  Lie algebras of type A, arXiv:1202.4227v1.


\end{thebibliography}
\end{document}